# Optimal targeting of nonlinear chaotic systems using a novel evolutionary computing strategy

Yudong Wang [a], Xiaoyi Feng [b], Xin Lyu [b], Zhengyang Li [b], Bo Liu [b,*]

[a] *Science and Technology on Space Physics Laboratory, Beijing, China*

[b] *Academy of Mathematics and Systems Science, Chinese Academy of Sciences, Beijing, China*

**Abstract**

Control of chaotic systems to given targets is a subject of substantial and well-developed research issue in nonlinear science, which can be formulated as a class of multi-modal constrained numerical optimization problem with multi-dimensional decision variables. This investigation elucidates the feasibility of applying a novel population-based metaheuristics labelled here as Teaching-learning-based optimization to direct the orbits of discrete chaotic dynamical systems towards the desired target region. Several consecutive control steps of small bounded perturbations are made in the Teaching-learning-based optimization strategy to direct the chaotic series towards the optimal neighborhood of the desired target rapidly, where a conventional controller is effective for chaos control. Working with the dynamics of the well-known Hénon as well as Ushio discrete chaotic systems, we assess the effectiveness and efficiency of the Teaching-learning-based optimization based optimal control technique, meanwhile the impacts of the core parameters on performances are also discussed. Furthermore, possible engineering applications of directing chaotic orbits are discussed.

**Highlights**

► Control of chaotic systems is formulated as constrained optimization problem.

► Teaching-learning-based optimization is to direct the chaotic series.

► The efficacy of TLBO based optimal control technique have been demonstrated.

[*] Corresponding author. Tel: +86-10-62614473; *E-mail address*: liub01@mails.tsinghua.edu.cn or bliu@amss.ac.cn (B. Liu)

# 1. Introduction

Chaos is a kind of characteristic of non-linear systems, which is a bounded unstable dynamic behavior that exhibits sensitive dependence on initial conditions and includes infinite unstable periodic motions. Although it appears to be stochastic, it occurs in a deterministic nonlinear system under deterministic conditions. In recently years, growing interests from physics, chemistry, biology, electronics, controls and instrumentation have stimulated the studies of chaos so as to improve the industrial and manufacturing systems and processes which exhibit chaotic phenomena [1]. Control of chaotic systems is one of important and well-developed research issues in nonlinear science [2-4], and it is possible that a minute perturbation of the control parameter could redirect chaos towards the desired region and stabilize it [5]. In the last decade, chaotic control has been shown to cover a wide spectrum of real world applications in engineering [6-17].

Since the pioneering work of Hübler's on chaos control in 1989 [18], a wide variety of approaches have been proposed for the appropriate control of chaotic systems. Most of the control techniques have been based on the *OGY* method [5, 19]. OGY exploits the exponential sensitivity of chaotic systems by using minute perturbations to direct the system towards a desired target in a short time. By extending the work of Ott et al. [5], Grebogi and Lai [20] described a method that converts the motion on a chaotic attractor to a desired attracting time periodic motion by making only small time dependent perturbations of a control parameter. This allows for a more generic choice of the feedback matrix and implementation to higher-dimensional systems. *Optimal control theory based approach* is an alternative approach for the control of chaotic dynamic systems [21]. Paskota et al. [22] applied the optimal control theory to calculate an open-loop controller and direct the orbit of a chaotic system towards the neighborhood of the desired target. Abarbanel et al. [3] demonstrated the use of an explicit single-step control method for directing a nonlinear system to the target orbit and maintaining it there. On the other hand, a few research studies have dealt with the control and synchronization problems of chaotic systems using the *variable structure control* scheme. Yu [23] examined the variable structure control strategy for the control of chaos in dynamic systems. In their study, by switching between two configurations of a perturbed parameter in chaotic systems, sliding regions can be created in which the desired performance lies. The stabilization and tracking of a periodic signal of the Rössler system have also been studied. Wang and Su [24] proposed an adaptive complementary variable structure control for chaotic synchronization. Based on Lyapunov's stability theory and the Babalat's lemma the proposed controller has been shown to render the synchronous error to zero. Fuh and Tung [25] presented an effective approach for controlling chaos by using a

*differential geometric method* which transformed nonlinear dynamics into linear form algebraically, so that linear control techniques can be used. It has been shown that the proposed method is not only able to control chaotic motion to a steady state but also to any desired periodic orbit. *Linear feedback methods* [26, 27] and *nonlinear feedback control* [28-30] are feasible approaches to direct chaotic systems to a steady state. [31] addressed the control of discrete-time chaotic dynamic systems using conventional feedback control strategies. Further, Chen and Dong [32] presented how to use the canonical feedback controllers to control chaotic trajectory of a continuous-time nonlinear system in converging to its equilibrium points and, more significantly, to its multiperiodic orbits including unstable limit cycles. Chen and Han [28] stabilize the controlled system at origin and synchronize two Genesio systems by designing a nonlinear feedback controller, whose stability could be easily guaranteed by using Hurwitz stability analysis approach. Based on the *adaptive control approach* [33], Liao and Tsai [34] constructed an adaptive observer-based driven system to synchronize the drive system whose dynamics are subjected to the system's disturbances and/or some unknown parameters. By appropriately selecting the observer gains, the synchronization and stability of the overall systems can be guaranteed with the Lyapunov approach. Fotsin and Bowong [35] addressed the problem of control and synchronization of coupled second-order oscillators. Firstly, they designed feedback controller to stabilize the system at its equilibrium. Then an adaptive observer was designed to synchronize the states of the master and slave oscillators using a single scalar signal corresponding to an observable state variable of the driving oscillator. Aghababa and Hashtarkhani [36] addressed the issue of synchronizing two different uncertain chaotic systems with unknown and different bounds via adaptive control method. Nijmeijer and Mareels [37] reformulated the chaotic synchronization as an *observer design problem*. Yang and Chen [38] provided some new observer-based criteria for discrete-time generalized chaos synchronization. Bai and Lonngren [39] utilized *active control theory* to synchronize a coupled Lorenz system. Tang and Wang [40] proposed an adaptive active control approach to make the states of two identical Chua's systems with unknown constant parameters to be asymptotically synchronized. Lu and Zhang [41] proposed the *backstepping design technique* for controlling Chen's chaotic attractor based on parameters identification. Wu and Lu [42] first designed an observer to identify the unknown parameter of Lü system, then applied the backstepping approach to control the uncertain Lü system to bounded points. Xu and Teo [43] considered the asymptotical stabilization problem of discrete chaotic systems based on the *impulsive control scheme*. By means of the Lyapunov stability theory and algebraic inequality techniques, sufficient conditions for asymptotical stability of the impulsive controlled discrete systems were obtained. Based on the impulsive control approach, Kemih et al. [44] addressed the satellite attitude control problem subjected to deterministic external perturbations which induced chaotic motions. Theorems on the stability of impulsive control systems were developed to find the conditions under which the chaotic systems can be asymptotically controlled to the origin by using impulsive control. As for the fuzzy approach [45], Poursamad and Davaie-Markazi [46] presented a robust *adaptive fuzzy control algorithm* for controlling unknown chaotic systems. The fuzzy system is designed to mimic an ideal controller, based on sliding-mode control. The robust controller is designed to compensate for the difference between the fuzzy controller and the ideal controller. The adaptive laws are derived in the Lyapunov sense to guarantee the stability of the controlled system. In addition to the above methods, some researches applied *optimization based methods* to direct chaos to targeted regions. From the viewpoint of optimization, control of chaotic systems could be formulated as multi-modal constrained numerical optimization problems [47-49]. Genetic algorithm [50], simplex-annealing strategy [51], Particle swarm optimization [52], and Differential Evolution [53] have been considered. Wang et al. [51] proposed an effective hybrid optimization strategy by combining the probabilistic jump search of simulated annealing with the convex polyhedron-based geometry search of Nelder-Mead Simplex method. The hybrid optimization strategy was applied to direct orbits of chaotic systems to a desired target region and to synchronize the two chaotic systems. Simulations results obtained on Hénon Map demonstrated the effectiveness of their hybrid approach.

In the past two decades, population-based optimization has attracted great attention from both academia and industry in many fields not limited in system science [54-62]. Recently, a new population-based metaheuristics, labeled as the Teaching-learning-based optimization (TLBO), has been proposed [63-67] as an alternative to genetic algorithm (GA) [68], particle swarm optimization (PSO) [69, 70] and Differential Evolution (DE) [71] for continuous optimization problems. The TLBO is inspired by the process of the teaching process and learning process of students in a class. In TLBO, firstly a population of solutions which is composed of teacher and students is initialized randomly, in which the most knowledgeable individual with the best fitness value is generally regarded as the teacher, while the remaining individuals in population are considered as students. Then the population is evolved to find optimal solutions through *teaching phase* in which the teacher helps the students to improve their grades as well as the *learning phase* in which the students improve their grades through interactions among themselves. Compared with GA, PSO and DE, TLBO has some attractive characteristics. It uses simple differential operation between teacher and students to create new candidate solutions, as well as to guide the search toward the most promising region. The conventional TLBO only contains one adjustable controlling parameter which facilitates easy tuning and implementation, while in GA, PSO and DE more parameters need to be set in appropriate manner so as to guarantee the searching performance. Nowadays, TLBO has attracted much attention and wide applications in different fields since its birth in 2011 [72, 73]. Application areas cover dynamic economic emission dispatch [74], structural optimization [75], power system [76], heat exchangers [77, 78], thermoelectric cooler [79], chaotic time series prediction [80], planning and scheduling [81-85], bioinformatics [86] and engineering optimization problems [87-95] etc., which demonstrate the effectiveness and efficiency of the TLBO based algorithms.

To date, there has been a lack of research study on TLBO for chaos control. The objective of this investigation is explicitly set out to fulfill this role. In this study, the TLBO is applied to direct the orbits of chaotic dynamical systems, which could be formulated as multimodal numerical optimization problems with high dimensions. Simulations results based on Hénon Map and Ushio Map are then obtained to verify the effectiveness and efficiency of TLBO, and the effects of some core parameters are also investigated.

## 2. Problem formulation

Consider the following discrete chaotic dynamic system:

$$\mathbf{x}(k+1) = \mathbf{f}(\mathbf{x}(k)), \quad k = 1, 2, \ldots, N \quad (1)$$

where state $\mathbf{x}(k) \in R^n$, and $\mathbf{f}: R^n \to R^n$ is continuously differentiable.

To direct the system towards a desired target, often minute perturbation $\mathbf{u}(k) \in R^n$ is added to the chaotic system. The system depicted in Eq. (1) can then be reformulated as follows:

$$\mathbf{x}(k+1) = \mathbf{f}(\mathbf{x}(k)) + \mathbf{u}(k), k = 0, 1, \cdots, N-1 \quad (2)$$

where $\mathbf{x}_0 \in R^n$ is an initial state of the system, the value of small perturbation $\mathbf{u}(k)$ is bounded by $\|\mathbf{u}(k)\| \leq \mu$, and $\mu$ is a positive real constant.

The goal is to determine a trajectory of consecutive perturbations $\mathbf{u}(k)$, $k = 0, 1, \ldots, N-1$ in guiding the state of the subsequent $N$ steps, i.e., $\mathbf{x}(N)$ towards the $\varepsilon$-neighborhood of desired target $\mathbf{x}_t$, i.e., $\|\mathbf{x}(N) - \mathbf{x}_t\| < \varepsilon$, where a local controller or conventional control technique that is effective for chaos control lies.

Without loss of generality, assuming that perturbation $\mathbf{u}(k)$ is imposed only on the first component of $\mathbf{f}$, the optimal targeting of nonlinear chaotic systems can thus be casted as a multi-dimensional constrained numerical optimization problem of the following form:

$$\min_{u(k)} \quad \|\mathbf{x}(N) - \mathbf{x}_t\| \tag{3}$$

$$s.t. \quad \begin{cases} x_1(k+1) = f_1(\mathbf{x}(k)) + u(k) \\ x_i(k+1) = f_i(\mathbf{x}(k)), \quad i = 2, \cdots, n \end{cases} \tag{4}$$

$$|u(k)| \leq \mu \tag{5}$$

$$\mathbf{x}(0) = \mathbf{x}_0 \tag{6}$$

In the present study, the TLBO is considered for solving the above minimization problem which involves finding a series of suitable $u(k)$, $k = 0,1,...,N-1$ by minimizing objective function (3); Meanwhile, constraints arising from chaotic dynamics (4), imposed on amplitude of perturbation (5), and stemmed from the initial chaotic states (6) should not be violated in the TLBO search. In the experiments, the chaotic dynamical system $\mathbf{f}: R^n \to R^n$ is then instantiated with the well-established Hénon and Ushio discrete chaotic systems.

## 3. Teaching-learning-based optimization (TLBO)

In this section, the TLBO approach is described. In the TLBO system, a population of solutions corresponding to a group of learners is initialized randomly. The most knowledgeable individual, which is analogous to the elite solution with the best fitness value in the search, behaves as the teacher, while the remaining individuals in population are considered as the learners or students. Each dimension of an individual solution in the TLBO models the grade of a different subject as attained by a teacher or learner. The population is then evolved to locate optimal solutions through a *teaching phase* in which the teacher helps the students to improve their grades as well as a *learning phase* where students improve their grades through interactions among themselves.

The $i$-th individual in the $d$-dimensional search space at generation $t$ can be represented as $X_i(t) = [x_{i,1}, x_{i,2},...,x_{i,d}]$, ($i = 1,2,...,NP$, where $NP$ denotes the size of the population. Remark: $X_i(t)$ is the decision variables $u(k)$, $k = 0,1,...,N-1$ of Eq. (4) which are bounded according to Eq. (5)). As the teacher is considered the most knowledgeable person, the best member $X_{best}(t)$ of the current population as defined by the objective function or fitness value is considered as the teacher. In problem minimization, the solution or individual with the smallest objective function value is thus regarded as the best member. At each generation *t*, the *teaching* and *learning* operations are applied on the learners, and a new population arises. Then, *comparison* takes place, and the corresponding individuals from both populations compete to comprise the next generation.

For each learner $X_i(t)$, according to the *teaching* operation, an updated learner $V_i(t) = [v_{i,1},...,v_{i,j},...,v_{i,d}]$ is generated by adding the weighted difference between the teacher and mean grade of learners to itself, which takes the following form:

$$V_i(t) = R \cdot * (X_{best}(t) - T \cdot * M(t)) + X_i(t) \qquad (7)$$

where the arithmetic operator $\cdot *$ denotes element-by-element multiplication. $R = [rand_1,...,rand_j,...,rand_d]$ is *d*-dimensional random weight vector which controls amplification of the differential variation $X_{best}(t) - T \cdot * M(t)$, and each element $rand_j$ is the *j*-th independent random number which is uniformly distributed in the range of [0, 1]. As previously described, $X_{best}(t)$, the base vector to model after, is the best member of the current population so that the finest traits of the teacher can be passed to the learners. $M(t) = [m_1,...,m_j,...,m_d]$ denotes the mean grade of the learners for each subject. $T$, known as the teaching factor which represents the aptitude of the teacher, is a *d*-dimensional random weight vector that controls the changes to the mean grades of learners. The value for each element of $T$ is then either 1 or 2, as recommended in Rao et al. [64].

After all the learners have completed the *teaching* phase, the *one to one selection* operator is then applied on each individual to decide whether the updated learner $V_i(t)$ or the original $X_i(t)$ would become a member of the population that would subsequently undergo the *learning* phase. Thus, for each target individual, a *new trial vector* $U_i(t) = [u_{i,1}(t),...,u_{i,d}(t)]$ is generated and assigned to the value $X_i(t)$ if the target learner could not improve itself in the teaching process; otherwise $U_i(t)$ is set to be $V_i(t)$.

For each *trial vector* of the *i-th* learner $U_i(t)$, through the *learning* phase, learners improve themselves by learning from others in the group, which is described by the following equation (for minimization problem):

$$W_i(t) = \begin{cases} U_i(t) + Rand.*(U_i(t) - U_j(t)), & if\ F(U_i(t)) < F(U_j(t)), \\ U_i(t) + Rand.*(U_j(t) - U_i(t)), & otherwise. \end{cases} \quad (8)$$

where $W_i(t)$ is the trial vector after learning for target individual $i$. *Rand* is a *d*-dimensional random weight vector which controls amplification of the differential variation. The subscript *j* of $U_j(t)$ denotes a randomly selected target individual *j* form the population $\{1,2,...,NP\}$ and also different from the current index *i*. $F(\ )$ means the objective evaluation function, which is the value of Eq.(3) in our study.

Finally, the *selection* arises to decide whether the *trial vector* $W_i(t)$ would be a member of the population of the next generation $t+1$. For a minimization problem, $W_i(t)$ is compared to $U_i(t)$ using the following one to one greedy based selection criterion:

$$X_i(t+1) = \begin{cases} W_i(t), & if\ F(W_i(t)) < F(U_i(t)), \\ U_i(t), & otherwise. \end{cases} \quad (9)$$

where $X_i(t+1)$ is the individual of the new population. The best individual of the current new population with the best objective value is then determined. If the objective value is better than the objective value of $X_{best}$, $X_{best}$ is then updated as the new best individual. The above operation is iterated until the stopping criterion is met, and $X_{best}$ is then the converged solution obtained. The procedure described above is considered as the standard version of TLBO. The key parameter in TLBO is *NP* (size of population). As described previously, the conventional TLBO only contains one adjustable controlling parameter which facilitates easy tuning and eases of implementation. This is in contrast to the GA, PSO and DE, which has more parameters that need to be appropriately defined in order to assert good and robust search performances. In our study, we shall provide a sensitivity analysis on the effects of the population size parameter on the search performances, specifically for chaotic control.

As for the problem to direct chaotic orbit, TLBO determines the consecutive control steps $(u(0), u(1), \cdots, u(N-1))$ which are bounded by the allowable perturbation amplitude while minimizing objective function (3). The control steps *N* denotes the dimension *d* of solution $X_i(t) = [x_{i,1}, x_{i,2}, ..., x_{i,d}]$, and $X_{best} = [u(0), u(1), ..., u(N-1)]$. In next section, the performance of TLBO is investigated on the direct chaotic orbit problem.

## 4. Numerical simulation on directing chaotic orbits

This section investigates the feasibility of applying the population-based metaheuristics TLBO for directing the orbits of discrete chaotic dynamical systems towards the desirable target region. Working with the dynamics of the well-established Hénon and Ushio discrete chaotic systems, we apply the TLBO to direct the orbits of discrete chaotic dynamical systems towards the desired target region; meanwhile the impacts of some parameters on controlling performance are also investigated. Furthermore, possible engineering applications of directing chaotic orbits are discussed.

### *4.1 Directing Hénon chaotic orbits*

As one of the most studied examples of the dynamical discrete-time chaotic system, Hénon Map is employed as the illustrating example in the present study, which can be described as follows:

$$\begin{cases} x_1(k+1) = -px_1^2(k) + x_2(k) + 1 \\ x_2(k+1) = qx_1(k) \end{cases} \tag{10}$$

The map has two parameters, $p$ and $q$, and for the classical Hénon map have values $p=1.4$, and $q=0.3$. On these classical values, the Hénon map is chaotic.

Target $\mathbf{x}_t$ is set to be a fixed point of the system (0.63135, 0.18941). Let $\mathbf{x}_0=(0, 0)$, and $u(k)$ is only appended to $x_1$ with bound $\mu=0.01$. The population size of TLBO is 50 and the maximum generation of TLBO is 1,000 (stopping condition), and the objective function is defined as $\|\mathbf{x}(N)-\mathbf{x}_t\|$. For different values of $N$, Table 1 reports the best, worst and the mean objective values as well as the standard deviations obtained from the 100 independent runs performed.

**(Insert Table 1 here)**

It could be concluded that the directing error (i.e. distance between the target and the ending state) in terms of 'Best objective value' decreases as $N$ increases, which suggests that a small number of consecutive control steps is insufficient to direct the Hénon chaotic trajectory towards the targeted region. Further, it can be observed that a consecutive control step of 9 give superior performance in terms of best, worst, mean objective values and standard deviations. In contrast, without above directing it would take a very long time to get close to the target starting from the same initial state, and the iteration (i.e. needed $N$) increases greatly as $\varepsilon$ decreases, see Table 2.

**(Insert Table 2 here)**

To test the performance of TLBO, TLBO is compared with standard Estimation of Distribution Algorithm (EDA) [96] and Particle Swarm Optimization (PSO) [52]. We adopt the simulation results from the literature. From Table 3, it can be seen that TLBO is superior to EDA and PSO in term of searching quality and derivation of the results. It could be concluded that TLBO is more effective and more robust on initial conditions.

**(Insert Table 3 here)**

Next, we investigate the impact of parameters $N$, $\mu$ and $\varepsilon$ on the directing errors of the Hénon chaotic orbits. The TLBO search is conducted 100 times independently for each combination of parameters, with the maximum number of generations set as 1000. If the resultant objective value in a run satisfies precision value $\varepsilon$, a success run is registered. Further, we define two performance measures labelled here as "succeed ratio (SR)" and "average valid evaluation number (AVEN)" which takes the following forms:

$$SR = \frac{N_v}{100} \times 100\% \quad (11)$$

$$AVEN = \frac{\sum_{i=1}^{N_v} n_i}{N_v} \quad (12)$$

where $N_v$ denotes the number of success runs among 100 independent runs, $n_i$ denotes the generation number of the $i$-th success run. Table 4 then reports the SR and AVEN results of the simulation runs obtained for the different combinations of $N$, $\mu$ and $\varepsilon$.

**(Insert Table 4 here)**

It could be seen from Table 4 that when a large directing error of precision $\varepsilon$ is considered, for instance $\varepsilon = 0.02$, TLBO is able to direct the chaotic orbit to the target region rapidly and is not affected by the configurations of the control step and perturbation bound parameters. However, as $\varepsilon$ is tightened (i.e., a higher directing error of precision is required, e.g., $\varepsilon = 0.001$ or $\varepsilon = 0.0001$), it can be seen that SR decreases while AVEN increases, respectively, on all combinations of $N$ and $\mu$. Fortunately, it is often unnecessary to use a small threshold value for local controller to be effective. Further, TLBO is noted to be effective where large control steps such as $N = 8$, $N = 9$, and $N = 10$, are used.

Next, we investigate the effect of perturbation bound parameter would affect the directing accuracy at each control step. Figs. 1-4 depict the mean objective values of directing errors over 100 independent runs for different configurations of $N$ and $\mu$. From these figures, the rapid decreasing trends of the mean objective values or errors are relatively steep and convergence happened at high precision. In addition, Figure 1 ($N = 7$) and Figure 4 ($N = 10$) also exhibits similar trends where the curves have been noted to descend faster on larger

$\mu$ values. On the contrary, Figure 2 ($N=8$) and Figure 3 ($N=9$) show opposite trends from Figures 1 and 4 where the curves plotted exhibits faster descend for small $\mu$ values. These seemingly paradoxical phenomena cannot be easily explained. A thoroughgoing landscape analysis of this problem may be necessary but is beyond the scope of the present study.

**(Insert Figure 1 here)**

**(Insert Figure 2 here)**

**(Insert Figure 3 here)**

**(Insert Figure 4 here)**

### *4.2 Directing Ushio chaotic orbits*

In what follows, we consider the Ushio discrete-time chaotic system [43, 97] which has the following form:

$$\begin{cases} x_1(k+1) = \alpha x_1(k) - x_1^3(k) + x_2(k) \\ x_2(k+1) = \beta x_1(k) \end{cases} \quad (13)$$

where $\alpha = 1.9$, and $\beta = 0.5$.

The target $\mathbf{x}_t$ is set to be (0, 0), the initial value $\mathbf{x}_0 = [0.6, -0.3]$, and $u(k)$ is only appended to $x_1$. Table 5 reports the best, worst and the mean objective values as well as the standard deviation of 100 independent runs for each combination of the control step and perturbation bound parameters. In this case, we consider the control steps from 6 to 20, and the perturbation bounds from 0.01 to 0.1.

**(Insert Table 5 here)**

From Table 5, first of all it can be observed that TLBO has been able to direct the trajectory of Ushio chaotic dynamic to the target region rapidly at high precision in nearly all cases, except the cases when both the control steps and perturbation bounds are configured as too small. The directing error is relatively larger when $N=6$ and perturbation bounds are 0.01 to 0.03. However, if perturbation bounds are extended to 0.04, then TLBO performs reasonably for $N=6$. Secondly, TLBO is shown to perform robustly with large control steps, which suggests it can handle high-dimension optimization problem. From this study, we conclude that TLBO is effectively and efficiently for directing chaotic orbits.

### *4.3 Possible applications of directing chaotic orbits*

Recently, control of chaotic systems to given targets has received increasing interests from both fundamental and applied researches. To improve the engineering systems and processes which exhibit chaotic phenomena, growing interests from physics, chemistry, biology, healthcare, electronics, and control have stimulated the studies on chaotic control [7, 98].

The motivations for controlling chaos in engineering practice are twofold: firstly, when chaotic dynamics are useful, the very nature of chaos, e.g., sensitive dependence on initial conditions, infinite unstable periodic motions, and ergodicity should be utilized instead of suppression. For instance, in chemical engineering, for the crucial equipment CSTR (Continuous Stirred Tank Reactor), by imposing chaotic dynamics on both the magnitude and the rate of change of the controlled variables, better mixing of reactants and full reaction could be well achieved [99].

Secondly, when chaotic dynamics or behaviors are undesired, the chaos should be redirected towards the desired region and be stabilized (i.e., decreased as much as possible) by a minute perturbation of the control parameter. Our study falls into this category in which the chaotic behaviors should be suppressed totally if possible. Taking the irregular vibrations in mechanical systems as the first example, the unfavorable vibrations arising from unbalanced rotators should be reduced as much as possible so as to avoid fatigue failures. Considering the fact that the motion of irregular vibrations is known to be chaotic, the control aim is then suppression of these undesirable vibrations by guiding or targeting the chaotic behaviors to a steady state [7, 25].

Anther representative example belonging to this category could be the control of chaotic attitude motion of the spacecraft with perturbation. Given the combination of forces from geomagnetic field, gravitational field, sunlight flux, etc., nearly all of the satellites in the spinning mode would exhibit chaotic motion [100-102]. Generally, the attitude motion of the spinning satellites could be characterized as nonlinear system of motion of rigid bodies. The attitude motion system is known to be several well-known chaotic dynamics, e.g., Lorenz equations, if periodically varying torque and long-term disturbing torque are selected. The objective of control is to redirect chaos towards the desired region (i.e., steady state) and stabilize, so that the chaotic behaviors in the spinning motion could be suppressed. Currently, a few approaches have been proposed for the appropriate control of chaotic attitude motion of the spacecraft with perturbation. For example, Mohammadbagheri and Yaghoobi [103] proposed a generalized predictive controller to suppress the chaos and regulate the state trajectory to desire fixed point. Kong et al. [100] proposed improved nonlinear relay control law based on position and velocities feedback to suppress the chaos in spacecraft attitude motion. Tsui and Jones [104] addressed the chaotic satellite attitude control problem by utilizing three techniques, i.e., delayed feedback control method, the Otani–Jones technique, and a higher dimensional variation of the OGY method, and their results showed that delayed feedback control method provided the most satisfactory control solution. Chen and Liu [105] applied the linear feedback method to control chaotic attitude motions of a magnetic rigid spacecraft with internal damping to the given fixed point. Abtahi et al. [106] investigated control of chaos for a Gyrostat satellite and designed OGY based method by using the linearization of the Poincaré map for suppression of chaos. Faramin and Ataei [107] investigated chaotic attitude maneuvers in a satellite for a range of parameters and designed back-stepping sliding mode method to ensure chaos suppression and achieve desired performance.

Based on the above investigation and analysis of the work on control of chaotic attitude motion of the spacecraft with perturbation, we would like to remark that a few control techniques have been investigated, e.g., the OGY method, optimal control theory based approach, nonlinear/linear feedback methods, adaptive control approach and backstepping design technique. However, in addition to the aforementioned methods, application of optimization based methods to control chaotic attitude motion of the spacecraft (i.e., direct chaos to targeted regions) is limited and relatively new, which is also our main motivation for writing this section on possible applications of directing chaotic orbit. And the proposed optimal targeting technique using the novel TLBO based evolutionary computing strategy is also necessary, and would be more useful for a practical user. It is also believed that the present study will provide benefits to a wide spectrum of real world applications of chaotic control in engineering, including maneuvering target tracking.

**5. Conclusion and further study**

Control of chaotic systems to given targets is the substantial research issue in both theory and application. By formulating it as a multi-dimensional constrained optimization problem, we investigated the feasibility of applying a novel population-based metaheuristics (Teaching-learning-based optimization, TLBO) for directing the orbits of discrete chaotic dynamical systems. Several consecutive control steps of small bounded perturbations were made in the Teaching-learning-based optimization strategy to direct the chaotic series towards the optimal neighborhood of the desired target rapidly, where a conventional controller was effective for chaos control. Working with the dynamics of the well-known Hénon as well as Ushio discrete chaotic systems, we assessed the effectiveness and efficiency of the Teaching-learning-based optimization based optimal control technique. First of all it could be observed that TLBO has been able to direct the trajectory of chaotic dynamic to the target region rapidly at high precision. Secondly, TLBO was shown to perform robustly with large control steps, which suggests it can handle high-dimension optimization problem. From this study, we concluded that TLBO was effectively and efficiently for directing chaotic orbits. Besides, TLBO was compared with standard Estimation of Distribution Algorithm (EDA) and Particle Swarm Optimization (PSO). It could be seen that TLBO was superior to EDA and PSO in term of searching quality and derivation of the results. Meanwhile, the impacts of the core parameters on control performance were also discussed. To the best of our knowledge, this is the first study to be reported on the use of TLBO for addressing problems involving chaos control.

In the past decades, Memetic Algorithms (MAs) have been demonstrated to converge to high-quality solutions more efficiently than their conventional counterparts on a wide range of benchmark problems as well as real-practice optimization problems [58, 108]. Among MAs, the Meta-Lamarckian learning [61] and the Probabilistic Memetic Framework [108] show outstanding searching performance on a variety of optimization problems. Besides, transfer learning concept provides a framework to utilize previously-acquired knowledge to solve new but similar problems much more quickly and effectively [109]. Our future work will investigate the performances of TLBO enhanced by Meta-Lamarckian learning scheme, the probabilistic memetic framework as well as transfer learning on the optimal targeting problems of nonlinear chaotic systems, respectively. As noticed that control of chaotic systems to given targets could be formulated as a class of multi-modal constrained numerical optimization problem with multi-dimensional decision variables. To further enhance the effectiveness and efficacy of the TLBO, a thoroughgoing landscape analysis of the formulated optimization problems may be necessary.

We hope the novel TLBO-based optimization methodology could serve as a vital control method for chaotic dynamic systems in more complex situations, e.g., control and synchronization of two different chaotic systems, control and synchronization of hyperchaotic systems, as well as targeting of nonlinear chaotic systems under uncertain parameters. It is also believed that the current study would provide benefits to a wide spectrum of real world applications of chaotic control in engineering practice, including control of chaotic attitude motion of the spacecraft with perturbation, maneuvering target tracking, and construction of secure communication system.

Our future work will investigate the performances of TLBO or improved TLBO on such substantial issues.

**Acknowledgements**

The authors wish to thank the anonymous referees and Editor-in-Chief Prof. Jie Lu for their constructive comments on earlier drafts of this paper. This work was supported in part by National Natural Science Foundation of China (Grant Nos. 71101139, 71390331, 11272062), National Science Fund for Distinguished Young Scholars of China (Grant No. 61525304). The authors are very grateful to Prof. Yew-Soon Ong (School of Computer Science & Engineering, Nanyang Technological University, Singapore) who has provided very useful comments and invaluable help to improve this study. The authors thank Royal Honored Prof. M.A. Keyzer (Faculty of Economics and Business Administration, SOW-VU, Vrije Universiteit Amsterdam, The Netherlands), Emeritus Prof. Yihui Jin, Prof. Ling Wang (Department of Automation, Tsinghua University), Prof. Jikun Huang (China Center for Agricultural Policy, Peking University) and Prof. Shouyang Wang (Academy of Mathematics and Systems Science, Chinese Academy of Sciences) for their support.

**Reference**

[1]   T. Kapitaniak, *Controlling chaos: theoretical and practical methods in non-linear dynamics*: Academic Press London, 1996.

[2]   T. Kapitaniak, "Continuous Control and Synchronization in Chaotic Systems," *Chaos Solitons & Fractals,* vol. 6, pp. 237-244, 1995.

[3]   H. D. I. Abarbanel, L. Korzinov, A. I. Mees, and I. M. Starobinets, "Optimal control of nonlinear systems to given orbits," *Systems & Control Letters,* vol. 31, pp. 263-276, Oct 10 1997.

[4]   F. O. Souza and R. M. Palhares, "Synchronisation of chaotic delayed artificial neural networks: an H-infinity control approach," *International Journal of Systems Science,* vol. 40, pp. 937-944, 2009.

[5]   E. Ott, C. Grebogi, and J. A. Yorke, "Controlling Chaos," *Physical Review Letters,* vol. 64, pp. 1196-1199, Mar 12 1990.

[6]   G. Chen and X. Dong, *From chaos to order : methodologies, perspectives, and applications*. Singapore ; River Edge, NJ: World Scientific, 1998.

[7]   A. L. Fradkov and R. J. Evans, "Control of chaos: Methods and applications in engineering," *Annual Reviews in Control,* vol. 29, pp. 33-56, 2005.

[8]   A. P. M. Tsui and A. J. Jones, "The control of higher dimensional chaos: comparative results for the chaotic satellite attitude control problem," *Physica D: Nonlinear Phenomena,* vol. 135, pp. 41-62, 2000.

[9]   T. Z. Rong and Z. Xiao, "Nonparametric interval prediction of chaotic time series and its application to climatic system," *International Journal of Systems Science,* vol. 44, pp. 1726-1732, Sep 1 2013.

[10]  C. F. Chuang, W. J. Wang, Y. J. Sun, and Y. J. Chen, "Convergence time guarantee for uncertain chaotic systems synchronisation," *International Journal of Systems Science,* vol. 44, pp. 1052-1060, Jun 1 2013.

[11]  W. L. Li, "Tracking control of chaotic coronary artery system," *International Journal of Systems Science,* vol. 43, pp. 21-30, 2012.

[12]  J. W. Liang, S. L. Chen, and C. M. Yen, "Identification and verification of chaotic dynamics in a missile

system from experimental time series," *International Journal of Systems Science,* vol. 44, pp. 700-713, Apr 1 2013.

[13] R. Raoufi and A. S. I. Zinober, "Smooth adaptive sliding mode observers in uncertain chaotic communication," *International Journal of Systems Science,* vol. 38, pp. 931-942, Nov 2007.

[14] C. F. Hsu, J. Z. Tsai, and C. J. Chiu, "Chaos synchronization of nonlinear gyros using self-learning PID control approach," *Applied Soft Computing,* vol. 12, pp. 430-439, Jan 2012.

[15] M. S. Tavazoei, M. Haeri, S. Jafari, S. Bolouki, and M. Siami, "Some Applications of Fractional Calculus in Suppression of Chaotic Oscillations," *IEEE Transactions on Industrial Electronics,* vol. 55, pp. 4094-4101, Nov 2008.

[16] G. Q. Li, P. F. Niu, W. P. Zhang, and Y. Zhang, "Control of discrete chaotic systems based on echo state network modeling with an adaptive noise canceler," *Knowledge-Based Systems,* vol. 35, pp. 35-40, Nov 2012.

[17] J. Park, J. Lee, and S. Won, "Adaptive sliding mode control for stabilization of stochastic chaotic systems with multi-constraints," presented at the 10th Asian Control Conference (ASCC), 2015.

[18] A. Hübler, "Adaptive control of chaotic systems," *Helv. Phys. Acta,* vol. 62, pp. 343-346, 1989.

[19] T. Shinbrot, E. Ott, C. Grebogi, and J. A. Yorke, "Using chaos to direct trajectories to targets," *Physical Review Letters,* vol. 65, p. 3215, 1990.

[20] C. Grebogi and Y. C. Lai, "Controlling chaotic dynamical systems," *Systems & Control Letters,* vol. 31, pp. 307-312, Oct 10 1997.

[21] G. Chen, "Optimal control of chaotic systems," *International Journal of Bifurcation and Chaos,* vol. 4, pp. 461-463, 1994.

[22] M. Paskota, A. Mees, and K. Teo, "Directing orbits of chaotic dynamical systems," *International Journal of Bifurcation and Chaos,* vol. 5, pp. 573-583, 1995.

[23] X. H. Yu, "Variable structure control approach for controlling chaos," *Chaos Solitons & Fractals,* vol. 8, pp. 1577-1586, Sep 1997.

[24] C. C. Wang and J. P. Su, "A new adaptive variable structure control for chaotic synchronization and secure communication," *Chaos Solitons & Fractals,* vol. 20, pp. 967-977, Jun 2004.

[25] C. C. Fuh and P. C. Tung, "Controlling Chaos Using Differential Geometric-Method," *Physical Review Letters,* vol. 75, pp. 2952-2955, Oct 16 1995.

[26] E. R. Hunt, "Stabilizing High-Period Orbits in a Chaotic System - the Diode Resonator," *Physical Review Letters,* vol. 67, pp. 1953-1955, Oct 7 1991.

[27] K. Pyragas and A. Tamasevicius, "Experimental Control of Chaos by Delayed Self-Controlling Feedback," *Physics Letters A,* vol. 180, pp. 99-102, Aug 30 1993.

[28] M. Y. Chen and Z. Z. Han, "Controlling and synchronizing chaotic Genesio system via nonlinear feedback control," *Chaos Solitons & Fractals,* vol. 17, pp. 709-716, Aug 2003.

[29] M. K. Ali and J. Q. Fang, "Synchronization of chaos and hyperchaos using linear and nonlinear feedback functions," *Physical Review E,* vol. 55, pp. 5285-5290, May 1997.

[30] X. X. Liao and G. R. Chen, "On feedback-controlled synchronization of chaotic systems," *International Journal of Systems Science,* vol. 34, pp. 453-461, Jun 10 2003.

[31] G. Chen and X. Dong, "On feedback control of chaotic nonlinear dynamic systems," *International Journal of Bifurcation and Chaos,* vol. 2, pp. 407-411, 1992.

[32] G. R. Chen and X. N. Dong, "On Feedback-Control of Chaotic Continuous-Time Systems," *IEEE Transactions on Circuits and Systems I-Fundamental Theory and Applications,* vol. 40, pp. 591-601, Sep 1993.

[33]   J. P. Su and C. C. Wang, "Synchronizing cascade-connected chaotic systems with uncertainties and breaking chaotic cryptosystems via a novel adaptive control scheme," *International Journal of Bifurcation and Chaos,* vol. 15, pp. 2457-2468, Aug 2005.

[34]   T. L. Liao and S. H. Tsai, "Adaptive synchronization of chaotic systems and its application to secure communications," *Chaos Solitons & Fractals,* vol. 11, pp. 1387-1396, Jul 2000.

[35]   H. Fotsin and S. Bowong, "Adaptive control and synchronization of chaotic systems consisting of Van der Pol oscillators coupled to linear oscillators," *Chaos Solitons & Fractals,* vol. 27, pp. 822-835, Feb 2006.

[36]   M. P. Aghababa and B. Hashtarkhani, "Synchronization of Unknown Uncertain Chaotic Systems Via Adaptive Control Method," *Journal of Computational and Nonlinear Dynamics,* vol. 10, Sep 2015.

[37]   H. Nijmeijer and I. M. Y. Mareels, "An observer looks at synchronization," *Circuits and Systems I: Fundamental Theory and Applications, IEEE Transactions on,* vol. 44, pp. 882-890, 1997.

[38]   X. S. Yang and G. R. Chen, "Some observer-based criteria for discrete-time generalized chaos synchronization," *Chaos Solitons & Fractals,* vol. 13, pp. 1303-1308, May 2002.

[39]   E. W. Bai and K. E. Lonngren, "Sequential synchronization of two Lorenz systems using active control," *Chaos Solitons & Fractals,* vol. 11, pp. 1041-1044, Jun 2000.

[40]   F. Tang and L. Wang, "An adaptive active control for the modified Chua's circuit," *Physics Letters A,* vol. 346, pp. 342-346, Oct 31 2005.

[41]   J. H. Lu and S. C. Zhang, "Controlling Chen's chaotic attractor using backstepping design based on parameters identification," *Physics Letters A,* vol. 286, pp. 148-152, Jul 23 2001.

[42]   X. Q. Wu and J. A. Lu, "Parameter identification and backstepping control of uncertain Lu system," *Chaos Solitons & Fractals,* vol. 18, pp. 721-729, Nov 2003.

[43]   H. L. Xu and K. L. Teo, "Stabilizability of discrete chaotic systems via unified impulsive control," *Physics Letters A,* vol. 374, pp. 235-240, Dec 28 2009.

[44]   K. Kemih, A. Kemiha, and M. Ghanes, "Chaotic attitude control of satellite using impulsive control," *Chaos, Solitons & Fractals,* vol. 42, pp. 735-744, 2009.

[45]   S. Y. Li, "Chaos control of new Mathieu-van der Pol systems by fuzzy logic constant controllers," *Applied Soft Computing,* vol. 11, pp. 4474-4487, Dec 2011.

[46]   A. Poursamad and A. H. Davaie-Markazi, "Robust adaptive fuzzy control of unknown chaotic systems," *Applied Soft Computing,* vol. 9, pp. 970-976, Jun 2009.

[47]   D. W. Wang and W. H. Ip, "Ant search based control optimisation strategy for a class of chaotic system," *International Journal of Systems Science,* vol. 36, pp. 951-959, Dec 15 2005.

[48]   L. D. S. Coelho and B. M. Herrera, "Fuzzy identification based on a chaotic particle swarm optimization approach applied to a nonlinear yo-yo motion system," *Ieee Transactions on Industrial Electronics,* vol. 54, pp. 3234-3245, 2007.

[49]   Q. L. Wei, D. R. Liu, and Y. C. Xu, "Policy iteration optimal tracking control for chaotic systems by using an adaptive dynamic programming approach," *Chinese Physics B,* vol. 24, Mar 2015.

[50]   X. Zhong, S. Shao, and J. Fang, "Directing orbits of chaotic dynamical systems using genetic algorithms (in Chinese)," *Control and Decision,* vol. 13, pp. 165-168 1998.

[51]   L. Wang, L. L. Li, and F. Tang, "Directing orbits of chaotic systems using a hybrid optimization strategy," *Physics Letters A,* vol. 324, pp. 22-25, Apr 5 2004.

[52]   B. Liu, L. Wang, F. Tang, and D. Huang, "Directing orbits of chaotic systems by particle swarm optimization," *Chaos Solitons & Fractals,* vol. 29, pp. 454-461, Jul 2006.

[53]   B. Liu, L. Wang, Y. H. Jin, D. X. Huang, and F. Tang, "Control and synchronization of chaotic systems by differential evolution algorithm," *Chaos Solitons & Fractals,* vol. 34, pp. 412-419, Oct 2007.

[54] X. S. Chen and Y. S. Ong, "A Conceptual Modeling of Meme Complexes in Stochastic Search," *IEEE Transactions on Systems Man and Cybernetics Part C-Applications and Reviews,* vol. 42, pp. 612-625, Sep 2012.

[55] M. N. Le, Y. S. Ong, Y. C. Jin, and B. Sendhoff, "A Unified Framework for Symbiosis of Evolutionary Mechanisms with Application to Water Clusters Potential Model Design," *IEEE Computational Intelligence Magazine,* vol. 7, pp. 20-35, Feb 2012.

[56] Y. Q. Wang, H. Z. Xie, X. Jiang, and B. Liu, "Intelligent Closed-Loop Insulin Delivery Systems for ICU Patients," *IEEE Journal of Biomedical and Health Informatics,* vol. 18, pp. 290-299, Jan 2014.

[57] B. Liu, L. Wang, Y. Liu, and S. Y. Wang, "A unified framework for population-based metaheuristics," *Annals of Operations Research,* vol. 186, pp. 231-262, Jun 2011.

[58] X. S. Chen, Y. S. Ong, M. H. Lim, and K. C. Tan, "A Multi-Facet Survey on Memetic Computation," *IEEE Transactions on Evolutionary Computation,* vol. 15, pp. 591-607, Oct 2011.

[59] Y. S. Ong, M. H. Lim, N. Zhu, and K. W. Wong, "Classification of adaptive memetic algorithms: a comparative study," *IEEE Transactions on Systems Man and Cybernetics Part B-Cybernetics,* vol. 36, pp. 141-152, 2006.

[60] Y. S. Ong, M. H. Lim, and X. S. Chen, "Memetic Computation-Past, Present & Future," *IEEE Computational Intelligence Magazine,* vol. 5, pp. 24-31, 2010.

[61] Y. S. Ong and A. J. Keane, "Meta-Lamarckian learning in memetic algorithms," *IEEE Transactions on Evolutionary Computation,* vol. 8, pp. 99-110, 2004.

[62] I. Giagkiozis, R. C. Purshouse, and P. J. Fleming, "An overview of population-based algorithms for multi-objective optimisation," *International Journal of Systems Science,* vol. 46, pp. 1572-1599, 2015/07/04 2015.

[63] R. V. Rao, V. J. Savsani, and D. P. Vakharia, "Teaching-learning-based optimization: A novel method for constrained mechanical design optimization problems," *Computer-Aided Design,* vol. 43, pp. 303-315, Mar 2011.

[64] R. V. Rao, V. J. Savsani, and J. Balic, "Teaching-learning-based optimization algorithm for unconstrained and constrained real-parameter optimization problems," *Engineering Optimization,* vol. 44, pp. 1447-1462, 2012.

[65] R. V. Rao and V. Patel, "Comparative performance of an elitist teaching-learning-based optimization algorithm for solving unconstrained optimization problems," *International Journal of Industrial Engineering Computations,* vol. 4, pp. 29-50, 2013.

[66] R. V. Rao and V. Patel, "An improved teaching-learning-based optimization algorithm for solving unconstrained optimization problems," *Scientia Iranica,* vol. 20, pp. 710-720, Jun 2013.

[67] R. Rao and V. Patel, "An elitist teaching-learning-based optimization algorithm for solving complex constrained optimization problems," *International Journal of Industrial Engineering Computations,* vol. 3, pp. 535-560, 2012.

[68] D. E. Goldberg, *Genetic algorithms in search, optimization, and machine learning*: Addison-Wesley, 1989.

[69] J. Kennedy, R. C. Eberhart, and Y. Shi, *Swarm Intelligence*. San Francisco: Morgan Kaufmann Publishers, 2001.

[70] L. Wang and B. Liu, *Particle swarm optimization and scheduling algorithms*, 2nd ed.: Beijing: Tsinghua University Press, 2011.

[71] R. Storn and K. Price, "Differential evolution - A simple and efficient heuristic for global optimization over continuous spaces," *Journal of Global Optimization,* vol. 11, pp. 341-359, 1997.

[72] G. Waghmare, "Comments on "A note on teaching-learning-based optimization algorithm"," *Information*

*Sciences,* vol. 229, pp. 159-169, Apr 20 2013.

[73] V. K. Patel and V. J. Savsani, "A multi-objective improved teaching–learning based optimization algorithm (MO-ITLBO)," *Information Sciences*.

[74] T. Niknam, F. Golestaneh, and M. S. Sadeghi, "theta-Multiobjective Teaching-Learning-Based Optimization for Dynamic Economic Emission Dispatch," *IEEE Systems Journal,* vol. 6, pp. 341-352, Jun 2012.

[75] T. Dede, "Optimum design of grillage structures to LRFD-AISC with teaching-learning based optimization," *Structural and Multidisciplinary Optimization,* vol. 48, pp. 955-964, Nov 2013.

[76] J. A. M. Garcia and A. J. G. Mena, "Optimal distributed generation location and size using a modified teaching-learning based optimization algorithm," *International Journal of Electrical Power & Energy Systems,* vol. 50, pp. 65-75, Sep 2013.

[77] R. V. Rao and V. Patel, "Multi-objective optimization of heat exchangers using a modified teaching-learning-based optimization algorithm," *Applied Mathematical Modelling,* vol. 37, pp. 1147-1162, Feb 1 2013.

[78] V. Patel and V. Savsani, "Optimization of a plate-fin heat exchanger design through an improved multi-objective teaching-learning based optimization (MO-ITLBO) algorithm," *Chemical Engineering Research and Design,* vol. 92, pp. 2371-2382, 11// 2014.

[79] R. Venkata Rao and V. Patel, "Multi-objective optimization of two stage thermoelectric cooler using a modified teaching–learning-based optimization algorithm," *Engineering Applications of Artificial Intelligence,* vol. 26, pp. 430-445, 1// 2013.

[80] R. G. Li, H. L. Zhang, W. H. Fan, and Y. Wang, "Hermite orthogonal basis neural network based on improved teaching-learning-based optimization algorithm for chaotic time series prediction," *Acta Physica Sinica,* vol. 64, Oct 20 2015.

[81] K. Xia, L. Gao, W. D. Li, and K. M. Chao, "Disassembly sequence planning using a Simplified Teaching-Learning-Based Optimization algorithm," *Advanced Engineering Informatics,* vol. 28, pp. 518-527, Oct 2014.

[82] H. Y. Zheng and L. Wang, "An effective teaching-learning-based optimisation algorithm for RCPSP with ordinal interval numbers," *International Journal of Production Research,* vol. 53, pp. 1777-1790, Mar 19 2015.

[83] Y. Xu, L. Wang, S. Y. Wang, and M. Liu, "An effective teaching-learning-based optimization algorithm for the flexible job-shop scheduling problem with fuzzy processing time," *Neurocomputing,* vol. 148, pp. 260-268, Jan 19 2015.

[84] J. Q. Li, Q. K. Pan, and K. Mao, "A discrete teaching-learning-based optimisation algorithm for realistic flowshop rescheduling problems," *Engineering Applications of Artificial Intelligence,* vol. 37, pp. 279-292, Jan 2015.

[85] T. Dokeroglu, "Hybrid teaching-learning-based optimization algorithms for the Quadratic Assignment Problem," *Computers & Industrial Engineering,* vol. 85, pp. 86-101, Jul 2015.

[86] D. L. Gonzalez-Alvarez, M. A. Vega-Rodriguez, and A. Rubio-Largo, "Finding Patterns in Protein Sequences by Using a Hybrid Multiobjective Teaching Learning Based Optimization Algorithm," *IEEE-ACM Transactions on Computational Biology and Bioinformatics,* vol. 12, pp. 656-666, May-Jun 2015.

[87] K. Yu, X. Wang, and Z. Wang, "An improved teaching-learning-based optimization algorithm for numerical and engineering optimization problems," *Journal of Intelligent Manufacturing,* pp. 1-13, 2014/05/13 2014.

[88] R. V. Rao and G. Waghmare, "Design optimization of robot grippers using teaching-learning-based optimization algorithm," *Advanced Robotics,* vol. 29, pp. 431-447, Mar 19 2015.

[89]  R. V. Rao, V. D. Kalyankar, and G. Waghmare, "Parameters optimization of selected casting processes using teaching-learning-based optimization algorithm," *Applied Mathematical Modelling,* vol. 38, pp. 5592-5608, Dec 1 2014.

[90]  B. D. Raja, R. L. Jhala, and V. Patel, "Multi-objective optimization of a rotary regenerator using tutorial training and self-learning inspired teaching-learning based optimization algorithm (TS-TLBO)," *Applied Thermal Engineering,* vol. 93, pp. 456-467, Jan 25 2016.

[91]  R. V. Rao and K. C. More, "Optimal design of the heat pipe using TLBO (teaching learning-based optimization) algorithm," *Energy,* vol. 80, pp. 535-544, Feb 1 2015.

[92]  W. W. Lin, D. Y. Yu, S. Wang, C. Y. Zhang, S. Q. Zhang, H. Y. Tian*, et al.*, "Multi-objective teaching-learning-based optimization algorithm for reducing carbon emissions and operation time in turning operations," *Engineering Optimization,* vol. 47, pp. 994-1007, Jul 3 2015.

[93]  Y. H. Cheng, "Estimation of Teaching-Learning-Based Optimization Primer Design Using Regression Analysis for Different Melting Temperature Calculations," *IEEE Transactions on Nanobioscience,* vol. 14, pp. 2-11, Jan 2015.

[94]  D. B. Chen, F. Zou, Z. Li, J. T. Wang, and S. W. Li, "An improved teaching-learning-based optimization algorithm for solving global optimization problem," *Information Sciences,* vol. 297, pp. 171-190, Mar 10 2015.

[95]  H. Boudjefdjouf, R. Mehasni, A. Orlandi, H. R. E. H. Bouchekara, F. De Paulis, and M. K. Smail, "Diagnosis of Multiple Wiring Faults Using Time-Domain Reflectometry and Teaching-Learning-Based Optimization," *Electromagnetics,* vol. 35, pp. 10-24, Jan 2 2015.

[96]  X. L. Huang, P. F. Jia, and B. Liu, "Controlling Chaos by An Improved Estimation of Distribution Algorithm," *Mathematical & Computational Applications,* vol. 15, pp. 866-871, 2010.

[97]  T. Ushio, "Chaotic Synchronization and Controlling Chaos Based on Contraction-Mappings," *Physics Letters A,* vol. 198, pp. 14-22, Feb 13 1995.

[98]  G. Chen and X. H. Yu, *Chaos control : theory and applications*. Berlin ; New York: Springer, 2003.

[99]  W. Wu, "Nonlinear bounded control of a nonisothermal CSTR," *Industrial and Engineering Chemistry Research,* vol. 39, pp. 3789–3798, 2000.

[100] L. Y. Kong, F. Q. Zhou, and J. Zou, "The control of chaotic attitude motion of a perturbed spacecraft," in *2006 Chinese Control Conference*, 2006, pp. 166-170.

[101] H. Wen, P. D. Jin, and Y. H. Hu, "Advances in dynamics and control of tethered satellite systems," *Acta Mechanica Sinica,* vol. 24, pp. 229-241, 2008.

[102] Z. Pang and D. Jin, "Experimental verification of chaotic control of an underactuated tethered satellite system," *Acta Astronautica,* vol. 120, pp. 287-294, 3// 2016.

[103] A. Mohammadbagheri and M. Yaghoobi, "Lorenz-Type Chaotic Attitude Control of Satellite through Predictive Control," in *2011 Third International Conference on Computational Intelligence, Modelling & Simulation*, 2011, pp. 147-152.

[104] A. P. M. Tsui and A. J. Jones, "The control of higher dimensional chaos: comparative results for the chaotic satellite attitude control problem," *Physica D: Nonlinear Phenomena,* vol. 135, pp. 41-62, 1/1/ 2000.

[105] L.-Q. Chen and Y.-Z. Liu, "Chaotic attitude motion of a magnetic rigid spacecraft and its control," *International Journal of Non-Linear Mechanics,* vol. 37, pp. 493-504, 4// 2002.

[106] S. M. Abtahi, S. H. Sadati, and H. Salarieh, "Nonlinear Analysis and Attitude Control of a Gyrostat Satellite with Chaotic Dynamics Using Discrete‐Time LQR‐OGY," *Asian Journal of Control,* 2016.

[107] M. Faramin and M. Ataei, "Chaotic attitude analysis of a satellite via Lyapunov exponents and its robust nonlinear control subject to disturbances and uncertainties," *Nonlinear Dynamics,* vol. 83, pp. 361-374,

2016.

[108] Q. H. Nguyen, Y. S. Ong, and M. H. Lim, "A Probabilistic Memetic Framework," *IEEE Transactions on Evolutionary Computation,* vol. 13, pp. 604-623, 2009.

[109] J. Lu, V. Behbood, P. Hao, H. Zuo, S. Xue, and G. Q. Zhang, "Transfer learning using computational intelligence: A survey," *Knowledge-Based Systems,* vol. 80, pp. 14-23, May 2015.

Table 1. Statistics performance of TLBO under different $N$ when $\mu = 0.01$ for Hénon chaotic orbits

| $N$ | Best objective value | Worst objective value | Mean objective value | Standard deviation |
|---|---|---|---|---|
| 6 | 0.08509 | 0.09331 | 0.09003 | 0.00163 |
| 7 | 0.00835 | 0.01287 | 0.01114 | 0.0008 |
| 8 | 0.00002 | 0.00052 | 0.00021 | 0.00009 |
| 9 | 0.00001 | 0.00036 | 0.00012 | 0.00007 |
| 10 | 0.00002 | 0.09152 | 0.00163 | 0.01042 |

Table 2. Iteration number needed without directing for Hénon chaotic orbits

| $\varepsilon$ | 0.02 | 0.001 | 0.00001 |
|---|---|---|---|
| Needed $N$ | 1188 | 17342 | 3356954 |

Table 3. Comparisons between TLBO and state-of-the-art algorithms under different $N$ when $\mu = 0.01$ for Hénon chaotic orbits

| $N$ | | TLBO | EDA [96] | PSO [52] |
|---|---|---|---|---|
| 6 | Mean | 0.09003 | 0.09784 | 0.09390 |
|   | Best | 0.08509 | 0.09405 | 0.09390 |
|   | Standard deviation | 0.00163 | 0.0028 | N/A |
| 7 | Mean | 0.01114 | 0.02091 | 0.01340 |
|   | Best | 0.00835 | 0.01385 | 0.01290 |
|   | Standard deviation | 0.0008 | 0.00766 | N/A |
| 8 | Mean | 0.00021 | 0.02838 | 0.00085 |
|   | Best | 0.00002 | 0.00124 | 0.00047 |
|   | Standard deviation | 0.00009 | 0.03057 | N/A |
| 9 | Mean | 0.00012 | 0.05913 | 0.00061 |
|   | Best | 0.00001 | 0.00078 | 0.00000 |
|   | Standard deviation | 0.00007 | 0.05706 | N/A |
| 10 | Mean | 0.00163 | 0.06831 | 0.01820 |
|    | Best | 0.00002 | 0.00138 | 0.00000 |
|    | Standard deviation | 0.01042 | 0.03425 | N/A |

Table 4 Average generation number and success ratio under different parameters for Hénon chaotic orbits

| $N$ | $\mu$ | $\varepsilon = 0.02$ | | $\varepsilon = 0.001$ | | $\varepsilon = 0.0001$ | |
|---|---|---|---|---|---|---|---|
| | | SR(%) | AVEN | SR(%) | AVEN | SR(%) | AVEN |
| $N = 7$ | 0.01 | 100 | 1.57 | 0 | - | 0 | - |
| | 0.02 | 100 | 1.43 | 0 | - | 0 | - |
| | 0.03 | 100 | 1.64 | 1 | 52 | 0 | - |
| $N=8$ | 0.01 | 100 | 1.81 | 100 | 27.49 | 10 | 494.60 |
| | 0.02 | 100 | 1.57 | 100 | 18.80 | 54 | 440.03 |
| | 0.03 | 100 | 1.64 | 97 | 54.67 | 18 | 438.22 |
| $N=9$ | 0.01 | 100 | 2.12 | 100 | 26.61 | 42 | 511.73 |
| | 0.02 | 100 | 2.14 | 100 | 49.75 | 22 | 427.09 |
| | 0.03 | 100 | 2.21 | 100 | 90.94 | 17 | 459.11 |
| $N=10$ | 0.01 | 98 | 30.9 | 97 | 71.77 | 30 | 491.03 |
| | 0.02 | 100 | 5.09 | 82 | 114.29 | 11 | 430.09 |
| | 0.03 | 100 | 3.29 | 95 | 235.53 | 5 | 329.20 |

Table 5. Statistics performance of TLBO under different $N$ and $\mu$ for Ushio chaotic orbits

| $N$ | $\mu$ | Best objective value | Mean objective value | Worst objective value | Standard deviation |
|---|---|---|---|---|---|
| 6 | 0.01 | 0.22967 | 0.29073 | 0.33347 | 0.02402 |
| | 0.02 | 0.19066 | 0.2026 | 0.20821 | 0.00423 |
| | 0.03 | 0.17504 | 0.19154 | 0.20136 | 0.0064 |
| | 0.04 | 0.00529 | 0.17614 | 0.19147 | 0.03261 |
| | 0.05 | 0.00735 | 0.14905 | 0.18728 | 0.05099 |
| | 0.06 | 0.00445 | 0.12614 | 0.1863 | 0.05629 |
| | 0.07 | 0.00429 | 0.09344 | 0.17197 | 0.05294 |
| | 0.08 | 0.00635 | 0.09035 | 0.17632 | 0.04865 |
| | 0.09 | 0.00835 | 0.0734 | 0.15915 | 0.03658 |
| | 0.10 | 0.00328 | 0.05541 | 0.13686 | 0.03918 |
| 7 | 0.01 | 0.0521 | 0.39847 | 0.57298 | 0.1058 |
| | 0.02 | 0.03381 | 0.04032 | 0.04699 | 0.00231 |
| | 0.03 | 0.0248 | 0.03376 | 0.04254 | 0.00288 |
| | 0.04 | 0.01295 | 0.02621 | 0.0378 | 0.00498 |
| | 0.05 | 0.01175 | 0.02166 | 0.0282 | 0.00378 |
| | 0.06 | 0.00521 | 0.017 | 0.0285 | 0.0052 |
| | 0.07 | 0.0058 | 0.01347 | 0.02062 | 0.00402 |
| | 0.08 | 0.00211 | 0.01132 | 0.01985 | 0.00402 |
| | 0.09 | 0.00185 | 0.00812 | 0.01552 | 0.00342 |
| | 0.10 | 0.00073 | 0.00526 | 0.01166 | 0.00267 |
| 8 | 0.01 | 0.2153 | 0.67065 | 1.05228 | 0.21688 |
| | 0.02 | 0.00171 | 0.00569 | 0.05566 | 0.00754 |
| | 0.03 | 0.00006 | 0.00137 | 0.00398 | 0.00072 |
| | 0.04 | 0.00004 | 0.00037 | 0.00157 | 0.00026 |
| | 0.05 | 0.00003 | 0.00036 | 0.00087 | 0.00019 |
|    | 0.06 | 0.00007 | 0.00037 | 0.00093 | 0.0002  |
|    | 0.07 | 0.00008 | 0.00038 | 0.00114 | 0.00019 |
|    | 0.08 | 0.00002 | 0.00054 | 0.00162 | 0.00034 |
|    | 0.09 | 0.00014 | 0.00074 | 0.00255 | 0.00051 |
|    | 0.10 | 0.00009 | 0.00082 | 0.00179 | 0.0004  |
| 9  | 0.01 | 0.01656 | 0.77553 | 1.21201 | 0.36457 |
|    | 0.02 | 0.00015 | 0.01542 | 0.17558 | 0.03613 |
|    | 0.03 | 0.0001  | 0.00067 | 0.00185 | 0.00036 |
|    | 0.04 | 0.00019 | 0.00073 | 0.00179 | 0.00038 |
|    | 0.05 | 0.00011 | 0.00069 | 0.00168 | 0.00037 |
|    | 0.06 | 0.00012 | 0.00071 | 0.00154 | 0.00033 |
|    | 0.07 | 0.00008 | 0.00097 | 0.00318 | 0.00059 |
|    | 0.08 | 0.00014 | 0.00105 | 0.00271 | 0.00065 |
|    | 0.09 | 0.00017 | 0.00113 | 0.00236 | 0.00056 |
|    | 0.10 | 0.0002  | 0.00166 | 0.01289 | 0.0018  |
| 10 | 0.01 | 0.09116 | 0.92267 | 1.15636 | 0.27832 |
|    | 0.02 | 0.00059 | 0.0064  | 0.12064 | 0.01694 |
|    | 0.03 | 0.0001  | 0.00212 | 0.02009 | 0.00372 |
|    | 0.04 | 0.00009 | 0.0013  | 0.00473 | 0.00087 |
|    | 0.05 | 0.00016 | 0.00122 | 0.00327 | 0.00068 |
|    | 0.06 | 0.00009 | 0.00118 | 0.00277 | 0.00061 |
|    | 0.07 | 0.00014 | 0.00121 | 0.0029  | 0.00075 |
|    | 0.08 | 0.00014 | 0.00169 | 0.01283 | 0.00182 |
|    | 0.09 | 0.00025 | 0.00161 | 0.00814 | 0.00122 |
|    | 0.10 | 0.00021 | 0.00167 | 0.00336 | 0.00083 |
| 11 | 0.01 | 0.15884 | 0.78568 | 0.90108 | 0.17073 |
|    | 0.02 | 0.00238 | 0.06544 | 0.71442 | 0.13012 |
|    | 0.03 | 0.0002  | 0.00236 | 0.01573 | 0.0026  |
|    | 0.04 | 0.00017 | 0.00174 | 0.00588 | 0.00115 |
|    | 0.05 | 0.00031 | 0.00166 | 0.00397 | 0.00095 |
|    | 0.06 | 0.00044 | 0.00207 | 0.0048  | 0.00105 |
|    | 0.07 | 0.00014 | 0.00239 | 0.00741 | 0.00141 |
|    | 0.08 | 0.00033 | 0.00206 | 0.00419 | 0.00103 |
|    | 0.09 | 0.00015 | 0.00208 | 0.00725 | 0.00136 |
|    | 0.10 | 0.00022 | 0.00237 | 0.00804 | 0.00145 |
| 12 | 0.01 | 0.03563 | 0.71389 | 0.74352 | 0.12315 |
|    | 0.02 | 0.00158 | 0.04919 | 0.5285  | 0.11043 |
|    | 0.03 | 0.00011 | 0.0044  | 0.03587 | 0.00522 |
|    | 0.04 | 0.0005  | 0.00324 | 0.01489 | 0.00243 |
|    | 0.05 | 0.00048 | 0.00208 | 0.0079  | 0.00128 |
|    | 0.06 | 0.00008 | 0.00179 | 0.00493 | 0.00113 |
|    | 0.07 | 0.00026 | 0.00196 | 0.00491 | 0.00118 |
|    | 0.08 | 0.00025 | 0.00194 | 0.00395 | 0.0009  |
|    | 0.09 | 0.00026 | 0.00197 | 0.00497 | 0.00103 |
|    | 0.10 | 0.00038 | 0.00248 | 0.0057  | 0.00144 |
| 13 | 0.01 | 0.06301 | 0.21708 | 0.66242 | 0.06884 |
|    | 0.02 | 0.00227 | 0.0205  | 0.08783 | 0.01908 |
|    | 0.03 | 0.00034 | 0.00524 | 0.01475 | 0.00331 |
|    | 0.04 | 0.00022 | 0.00405 | 0.02393 | 0.00358 |
|    | 0.05 | 0.00054 | 0.0039  | 0.0084  | 0.00205 |
|    | 0.06 | 0.00029 | 0.0028  | 0.0087  | 0.00172 |
|    | 0.07 | 0.00039 | 0.0033  | 0.00743 | 0.00168 |
|    | 0.08 | 0.00057 | 0.00354 | 0.00971 | 0.00204 |
|    | 0.09 | 0.00038 | 0.00396 | 0.01303 | 0.00265 |
|    | 0.10 | 0.00045 | 0.00387 | 0.01522 | 0.00274 |
| 14 | 0.01 | 0.00685 | 0.06672 | 0.19632 | 0.03982 |
|    | 0.02 | 0.00425 | 0.02248 | 0.04316 | 0.01226 |
|    | 0.03 | 0.00106 | 0.00978 | 0.03742 | 0.00862 |
|    | 0.04 | 0.00026 | 0.00323 | 0.01149 | 0.00248 |
|    | 0.05 | 0.00026 | 0.00151 | 0.0056  | 0.00103 |
|    | 0.06 | 0.00022 | 0.00147 | 0.00534 | 0.00098 |
|    | 0.07 | 0.00024 | 0.00168 | 0.00978 | 0.00148 |
|    | 0.08 | 0.00015 | 0.00161 | 0.00432 | 0.00101 |
|    | 0.09 | 0.00032 | 0.00224 | 0.00666 | 0.00134 |
|    | 0.10 | 0.00055 | 0.00246 | 0.00607 | 0.00126 |
| 15 | 0.01 | 0.00868 | 0.02208 | 0.04581 | 0.00956 |
|    | 0.02 | 0.0007  | 0.00465 | 0.00727 | 0.00137 |
|    | 0.03 | 0.00082 | 0.00245 | 0.00383 | 0.00087 |
|    | 0.04 | 0.00009 | 0.00135 | 0.0039  | 0.00079 |
|    | 0.05 | 0.00006 | 0.001   | 0.00251 | 0.00067 |
|    | 0.06 | 0.00025 | 0.00115 | 0.00489 | 0.00082 |
|    | 0.07 | 0.0001  | 0.00127 | 0.0051  | 0.00087 |
|    | 0.08 | 0.00018 | 0.00134 | 0.00411 | 0.00094 |
|    | 0.09 | 0.00021 | 0.00156 | 0.00456 | 0.00095 |
|    | 0.10 | 0.00028 | 0.00201 | 0.00607 | 0.00135 |
| 16 | 0.01 | 0.00109 | 0.00652 | 0.00971 | 0.00215 |
|    | 0.02 | 0.00028 | 0.00205 | 0.00567 | 0.00149 |
|    | 0.03 | 0.00071 | 0.00161 | 0.00387 | 0.00076 |
|    | 0.04 | 0.00002 | 0.0012  | 0.00315 | 0.00075 |
|    | 0.05 | 0.00014 | 0.00119 | 0.00294 | 0.00062 |
|    | 0.06 | 0.00012 | 0.00152 | 0.00362 | 0.00081 |
|    | 0.07 | 0.00045 | 0.00141 | 0.00441 | 0.00092 |
|    | 0.08 | 0.00013 | 0.00189 | 0.00805 | 0.00138 |
|    | 0.09 | 0.00015 | 0.00177 | 0.00534 | 0.00097 |
|    | 0.10 | 0.00036 | 0.0021  | 0.00489 | 0.00119 |
| 17 | 0.01 | 0.00003 | 0.00053 | 0.00145 | 0.00031 |
|    | 0.02 | 0.00025 | 0.00124 | 0.00329 | 0.00074 |
|    | 0.03 | 0.00021 | 0.00124 | 0.00434 | 0.0008  |
| | 0.04 | 0.00002 | 0.00121 | 0.00303 | 0.0007 |
| | 0.05 | 0.00019 | 0.00132 | 0.00362 | 0.00075 |
| | 0.06 | 0.0003 | 0.00156 | 0.00393 | 0.00079 |
| | 0.07 | 0.00036 | 0.00142 | 0.00394 | 0.00075 |
| | 0.08 | 0.00055 | 0.00171 | 0.00403 | 0.00073 |
| | 0.09 | 0.00024 | 0.00212 | 0.00625 | 0.00132 |
| | 0.10 | 0.00041 | 0.00194 | 0.00438 | 0.00094 |
| 18 | 0.01 | 0.00007 | 0.00113 | 0.00346 | 0.00072 |
| | 0.02 | 0.00009 | 0.00098 | 0.00279 | 0.0006 |
| | 0.03 | 0.00028 | 0.00139 | 0.00265 | 0.00066 |
| | 0.04 | 0.0002 | 0.00142 | 0.00289 | 0.00071 |
| | 0.05 | 0.00037 | 0.00151 | 0.0032 | 0.0007 |
| | 0.06 | 0.0001 | 0.00156 | 0.00377 | 0.00096 |
| | 0.07 | 0.00031 | 0.00194 | 0.004 | 0.00095 |
| | 0.08 | 0.0003 | 0.00173 | 0.0035 | 0.00087 |
| | 0.09 | 0.00041 | 0.00214 | 0.00719 | 0.00137 |
| | 0.10 | 0.00045 | 0.00236 | 0.00656 | 0.00114 |
| 19 | 0.01 | 0.00014 | 0.00089 | 0.00255 | 0.00058 |
| | 0.02 | 0.00026 | 0.00123 | 0.00306 | 0.00061 |
| | 0.03 | 0.00015 | 0.00154 | 0.004 | 0.00085 |
| | 0.04 | 0.00021 | 0.00172 | 0.00419 | 0.00089 |
| | 0.05 | 0.00014 | 0.00177 | 0.00425 | 0.00101 |
| | 0.06 | 0.00029 | 0.00194 | 0.00446 | 0.00097 |
| | 0.07 | 0.00029 | 0.00187 | 0.00538 | 0.00106 |
| | 0.08 | 0.00041 | 0.00196 | 0.00531 | 0.00099 |
| | 0.09 | 0.00026 | 0.00255 | 0.007 | 0.00143 |
| | 0.10 | 0.00027 | 0.00221 | 0.00501 | 0.00118 |
| 20 | 0.01 | 0.00025 | 0.00147 | 0.00309 | 0.00074 |
| | 0.02 | 0.00035 | 0.00139 | 0.00318 | 0.00067 |
| | 0.03 | 0.00026 | 0.00168 | 0.00519 | 0.00114 |
| | 0.04 | 0.00011 | 0.00194 | 0.00391 | 0.00103 |
| | 0.05 | 0.00043 | 0.00201 | 0.00584 | 0.00122 |
| | 0.06 | 0.00038 | 0.00185 | 0.00409 | 0.00086 |
| | 0.07 | 0.00014 | 0.00208 | 0.00521 | 0.00116 |
| | 0.08 | 0.00046 | 0.00221 | 0.00469 | 0.00109 |
| | 0.09 | 0.00038 | 0.00217 | 0.00727 | 0.00133 |
| | 0.10 | 0.00032 | 0.00225 | 0.00687 | 0.00133 |

**Figure captions**

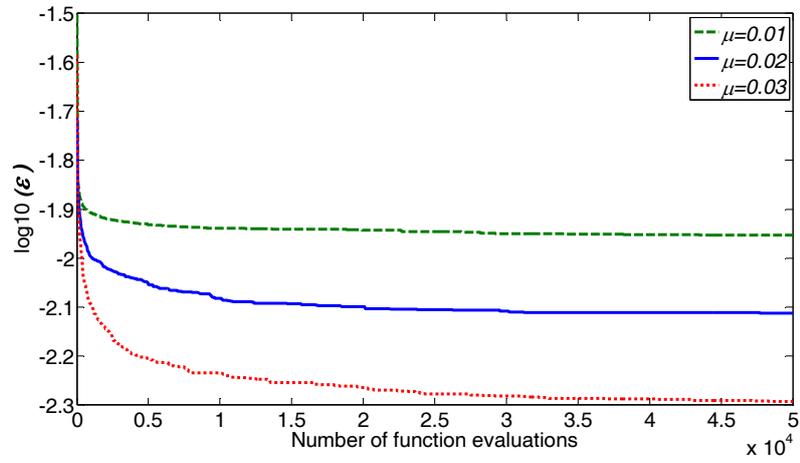

Figure 1 Descending curves of mean objective value when $N = 7$ for Hénon chaotic orbits

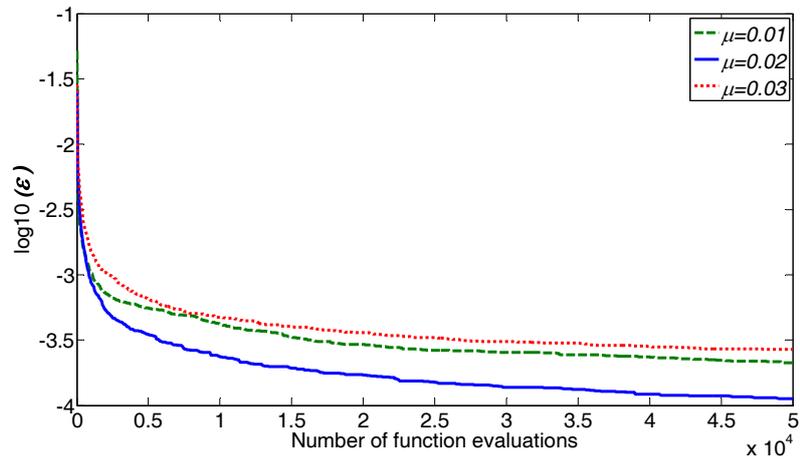

Figure 2 Descending curves of mean objective value when $N = 8$ for Hénon chaotic orbits

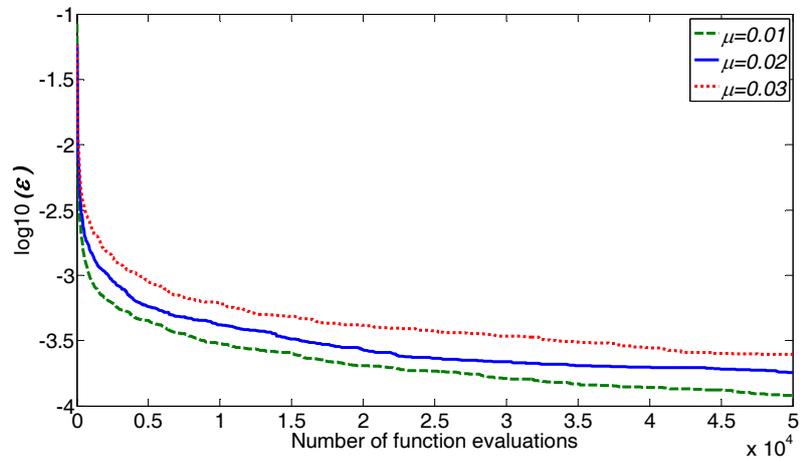

Figure 3 Descending curves of mean objective value when $N = 9$ for Hénon chaotic orbits

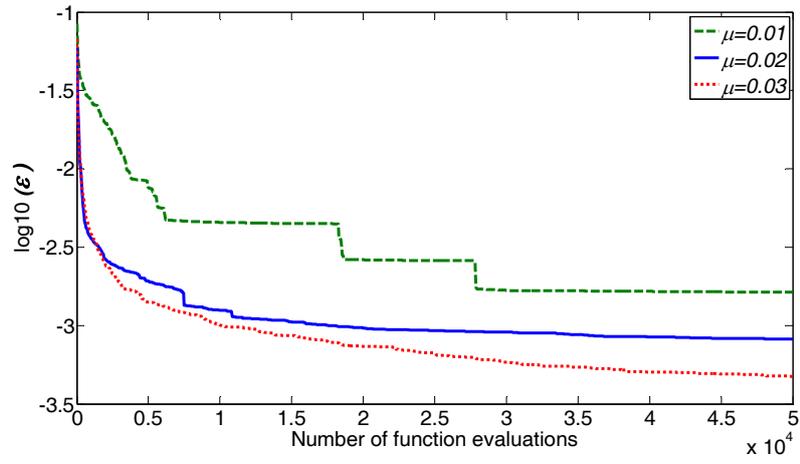

Figure 4 Descending curves of mean objective value when $N = 10$ for Hénon chaotic orbits